\documentclass[11pt]{article}

\newcommand{\ei}{E_{i}}
\newcommand{\eone}{E_{1}}
\newcommand{\etwo}{E_{2}}
\newcommand{\ethree}{E_{3}}
\newcommand{\tdx}{\textrm{dx}}
\newcommand{\tdy}{\textrm{dy}}
\newcommand{\ty}{\textrm{y}}
\newcommand{\tx}{\textrm{x}}

\newcommand{\elo}{\mathcal{L}^{1}}
\newcommand{\elt}{\mathcal{L}^{2}}
\newcommand{\ed}{E_{\Delta}}

\newcommand{\be}{\begin{eqnarray}}
\newcommand{\ee}{\end{eqnarray}}

\newcommand{\R}{\mathbb{R}}

\newcommand{\om}{\Omega}

\newcommand{\eps}{\epsilon}

\newcommand{\sca}{\mathcal{A}}

\newcommand{\scl}{\mathcal{L}}

\usepackage{amsfonts}
\usepackage{amsmath}
\usepackage{amsthm}
\usepackage{epsfig}
\usepackage{amssymb}
\usepackage{psfrag}
\newtheorem{thm}{Theorem}[section]
\newtheorem{prop}{Proposition}[section]
\newtheorem{lem}{Lemma}[section]

\newtheorem{defn}{Definition}[section]
\newtheorem{rem}[prop]{Remark}
\numberwithin{equation}{section}
\title{Local minimizers and low energy paths in a model of material microstructure with a surface energy term}
\author{J.J. Bevan\footnote{Department of Mathematics, University of Surrey, Guildford, Surrey,  GU2 7XH, UK.  email: j.bevan@surrey.ac.uk}}
\begin{document}
\markright{Local minimizers and low energy paths}
\maketitle
\begin{abstract}
 A family of integral functionals $\mathcal{F}$ which, in a simplified way,  model material microstructure occupying a two-dimensional domain $\om$ and which take account of surface energy and a variable well-depth is studied.   It is shown that there is a critical well-depth, whose scaling with the surface energy density and domain dimensions is given, below which the state $u=0$ is the global minimizer of a typical $F$ in $\mathcal{F}$.  It is also shown that $u=0$ is a strict local minimizer of $F$ in the sense that if $v \neq 0$ is admissible and either $||v||_{L^{2}(\om)}$ or $\elt(\{(x,y) \in \om: \ |v_{y}|(x,y) \geq 1\})$ is sufficiently small (with quantitative bounds given in terms of the parameters appearing in the energy functional $F$) then $F(v) > F(0)$.  Low energy paths between $u=0$ and the global minimizer (in the case of a sufficiently large well-depth) are given such that the cost of introducing sets $\{(x,y) \in \om: \ |v_{y}(x,y)| \geq 1\}$ of positive measure into the domain $\om$ may be made arbitrarily small.
\end{abstract}

\section{Introduction}
 

The energy functionals we shall consider in this paper are related to the one used in \cite{KM94} but with some important differences.  The original Kohn-M\"{u}ller functional is 
\[E_{KM}(u)= \int_{\om}\eps^{2} u_{yy}^{2}+ (u_{y}^{2}-1)^{2} +u_{x}^{2} \, \tdx,\]
where $\om=[0,L] \times [0,1]$, $L > 0$, $\eps >0$ is a small parameter which is sometimes referred to as the surface energy density, and $\tdx$ is shorthand for $d\elt(\textrm{x})$.

Provided suitable boundary conditions are imposed, the global minimizers of $E_{KM}$ model in a simplified way the fine scale microstructures that are observed to some degree at austenite-martensite interfaces in shape memory alloys.
The second derivative term tempers the oscillations in the $y-$direction that any globally minimizing sequence will develop.  Using a subtle argument Kohn and M\"{u}ller showed that the global minimization can be viewed as a straightforward competition between the term in $u_{x}^{2}$, effectively a measure of `spread', and a version of surface energy derived from the terms in $\int_{\om}u_{yy}^{2}\,\tdx$ and $\int_{\om}(u_{y}^{2}-1)^{2}\,\tdx$, among functions $u$ with $|u_{y}|=1$ a.e.  The result is that the infimum of the energy scales in $\eps$ as though it were evaluated at the now well-known branched microstructure (see \cite{KM94} for details).

However, the functional $E_{KM}$ is less useful in understanding the role of $u=0$ as a local minimizer and other details of the energy landscape.  To this end we introduce new functionals $\eone(\cdot; \eps, \Delta)$, $\etwo(\cdot; \eps, \Delta)$ and $\ethree(\cdot; \eps, \Delta)$ based on $E_{KM}$ but with some extra features.  The family of functionals $\mathcal{F}$ mentioned above consists of 
$\eone(\cdot; \eps, \Delta)$, $\etwo(\cdot; \eps, \Delta)$ and $\ethree(\cdot; \eps, \Delta)$ defined below in \eqref{E1}, \eqref{E2} and \eqref{E3} respectively as the parameters $\eps$ and $\Delta$ vary.
\begin{defn} \begin{eqnarray}
\label{E1} \eone(u; \eps, \Delta)& = & \int_{\om}\epsilon^{2}u_{yy}^{2}+u_{x}^{2}\,\emph{\tdx}+\Delta \elt(A(u)) \\
\label{E2}\etwo(u;\eps,\Delta) & = & \int_{\om}\epsilon^{2}|\nabla(u_{y})|^{2}+u_{x}^{2}\,\emph{\tdx}+\Delta \elt(A(u)) \\
\label{E3} \ethree(u;\eps,\Delta) & = &  \int_{\om} \eps^{2}|\nabla(\nabla(u))|^{2}+u_{x}^{2} \,\emph{\tdx} +\Delta \elt(A(u)) \\
\nonumber B(u) & = & \{(x,y) \in \om: \ |u_{y}(x,y)| \geq 1 \} \\
\nonumber  A(u) & = & \om \setminus B(u). \end{eqnarray}
\end{defn}

In each case the functionals take the form
\[E_{i}(u) = S_{i}(u) + \int_{\om}W_{\Delta}(\nabla u) \,\tdx,\]
where
\[W_{\Delta}(a,b)= a^{2}+ \Delta \chi_{(-1,1)}(b)\]
and where $S_{i}(u)$ is a surface energy term.    In particular, $S_{3}(u)$ will henceforth be written in the more conventional way
\[S_{3}(u) = \int_{\om}|D^{2}u|^{2}\,\tdx.\]

In the following we suppress the dependence of the $E_{i}$ on $\eps$ and $\Delta$ for brevity.
We study the behaviour of each $E_{i}$ in the class $\sca_{i}$ of admissible maps defined below.
First we define the subclass of $W_{p}^{1,2}(\om,\R) \subset W^{1,2}(\om,\R)$ of functions satisfying periodic boundary conditions (in the sense of trace) at the top and bottom of the domain $\om$ by
\begin{equation}\label{w12p} W_{p}^{1,2}(\om,\R)= \{u \in W^{1,2}(\om,\R): \ u(x,1)=u(x,0) \ \ 0 \leq x \leq 1 \}.\end{equation}
Then 
\begin{eqnarray}\label{a1}\sca_{1} &  = & \{u \in W_{p}^{1,2}(\om,\R): \, u_{yy} \in L^{2}(\om; \R),  u(0,y)=0 \  0 \leq y \leq 1\} \\
\label{a2} \sca_{2} & = & \{u \in W_{p}^{1,2}(\om,\R): \, u_{y} \in W^{1,2}(\om; \R),  u(0,y)=0 \  0 \leq y \leq 1\} \\
\label{a3} \sca_{3} & = & \{u \in W_{p}^{1,2}(\om,\R): \, \nabla u \in W^{1,2}(\om; \R),  u(0,y)=0 \,  0 \leq y \leq 1\}
\end{eqnarray}

The new features of these models relative to the original Kohn-M\"{u}ller functional are summarized here and discussed below.   They are:
\begin{itemize}
\item[(i)] a variable well-depth $\Delta$ (see below for its definition); 
\item[(ii)] a convex potential in a neighbourhood of $(0,0)$;
\item[(iii)] the possibility of a cost, which may be zero,  associated with the appearance of sets $B(u)$ of positive measure in $\om$, analagous to a lower bound on the cost of `nucleation' of martensite in austenite. 
\end{itemize}

The term $\elt(A(u))$ is intended to mimic the behaviour of the term $\int_{\om} (u_{y}^{2}-1)^{2} \tdx$ appearing in the Kohn-M\"{u}ller functional $E_{KM}$ in the following sense.  In order for $E_{i}(u)$ to approach its infimum $|u_{y}| < 1$ can only occur on a set of small measure: `most' values of $|u_{y}|$ will be near to or larger than 1.   Looking at $E_{KM}$ and referring to the argument they give we see that most values of the gradient $u_{y}$ of the global minimizer in that case will be near $\pm 1$.  The price we pay for replacing $\int_{\om} (u_{y}^{2}-1)^{2} \tdx$ with a term proportional to $\elt(A(u))$ is that large values of $|u_{y}|$ are not penalized as they would be in $E_{KM}$.  But this turns out not to matter a great deal, as will become clear later.  In fact, the nature of the global minimizer of $\eone$ in $\sca_{1}$ in a scaling sense can be deduced from the Kohn-M\"{u}ller argument when $\Delta$ is large enough, although we do not pursue this in the present work.  When $\Delta$ is in the range $(0,C\eps L^{-1})$ for an appropriate dimensionless constant $C$ it happens that $u=0$ is the global minimizer of $E_{i}$ in  $\sca_{i}$ for $i=1,2,3$.    This behaviour with respect to varying the well-depth $\Delta$ may simply be an artefact of the low dimensions and the choice of boundary conditions.  In any case it will be necessary to know just how large $\Delta$ needs to be before $u=0$ ceases to be the global minimizer of $E_{i}$ in $\sca_{i}$ for $i=1,2,3$.


The idea of introducing a well-depth 
\[\Delta =W_{\Delta}(0,0)-W_{\Delta}(0,\pm 1 )\] 
comes from the Ball-James theory of martensitic phase transformations.  The theory asserts that the stored-energy potential should change in a certain way as the temperature changes; we synthesize this by varying $\Delta$, with $\Delta=0$ corresponding to high temperature stored energy function and $\Delta=1$, say, to low temperature.  The global minimizer in the case $\Delta=0$ is $u=0$ which, in the full three-dimensional models, would be referred to as Austenite.   See \cite{BJ87}, \cite{BJ92} for further details.  We do not attempt to introduce dynamics.  It is shown in Section \ref{local} that $u=0$ is always a local minimizer of $\ei$ regardless of the size of $\Delta$; the only effect $\Delta$ has on local minimality is through the size of the neighbourhood $\mathcal{N}_{i}$, say,  of $u=0$ in $\sca_{i}$ on which $u=0$ satisfies $\ei(v) \geq \ei(0)$ for all $v$ in $\mathcal{N}_{i}$ for $i=1,2,3$.   See Section \ref{local}, and in particular Theorem \ref{t2},  for details.   When the well-depth is large enough it also appears in the scaling of the energy associated with the global minimizer; we do not pusue this in the present work.   We remark that in the case of $\eone$ one could follow the method given in \cite{KM94} with only minor changes.

Note that the potential   
\begin{equation*}W_{\Delta}(a,b)=a^{2}+\Delta \chi_{(-1,1)}(b)\end{equation*} 
is convex in a neighbourhood of $(a,b)=(0,0)$.  This is sufficient to establish that $u=0$ is an $L^{2}$-local minimizer of $\ei$ in $\sca_{i}$.  It is not necessary, though, as examples of Taheri show \cite[Section 4]{Ta00}.   In fact the potential $W_{\Delta}$ can be bounded below by a strongly convex potential in a neighbourhood of zero by `borrowing' some surface energy and applying a suitable Poincar\'{e} inequality.   For example, for any $\lambda \in (0,1)$,
\[\ed(u) \geq \int_{\om} (1-\lambda)\eps^{2} |\nabla u_{y}|^{2}  + C_{\om} \lambda \eps^{2}u_{y}^{2} +W_{\Delta}(\nabla u) \,\textrm{dx};\]
so that $\overline{W}_{\Delta}(s,t):=W_{\Delta}(s,t)+C_{\om} \lambda \eps^{2} t^{2}$ is strongly convex in a neighbourhood of $(0,0)$.  If we are allowed to vary the domain height, for example by taking 
\[\om_{h} = [0,L] \times [0,h]\]
one can use the same procedure to bound the Kohn-M\"{u}ller potential below by a potential that is strongly convex at $(0,0)$.  
The reason this works is that the Poincar\'{e} constant $C_{\om_{h}}$, say, becomes large when $h$ is small.  See Taheri \cite{Ta00} for other interesting examples. When $h=1$ the Poincar\'{e} constant is not large enough for this trick to work, from which it follows easily that $u=0$ is not a local minimizer of $E_{KM}$ in $\sca$ (provided $\eps$ is sufficiently small).
We note that the idea of combining the strong convexity of $W$ with higher order terms in order to guarantee local minimality has been studied before \cite{Ta00}, \cite{Ba05}.  Although we don't use these arguments directly they are, in view of the comments above, one of the main reasons that we can expect $u=0$ to be a local minimizer of $\ei$ in  $\sca_{i}$ for $i=1,2,3$.

The results of Section \ref{local} are based on an apparently new inequality which relates all three terms appearing in $\eone$ and which, together with more standard estimates, yields explicit functions $r(\eps,\Delta)$ and $s(\eps,\Delta)$ (given in \eqref{r} and \eqref{s} respectively) such that 
\begin{equation}\label{rs} ||u||_{L^{2}(\om)} < r(\eps,\Delta) \  \textrm{or} \ \scl^{2}(B(u)) <  s(\eps,\Delta) \implies \ed(u)  > \ed(0), \end{equation}
provided $u \neq 0$.  It is doubtful whether the scalings involved are optimal for reasons explained in Section \ref{local}.  Even so, it is still a stronger and more explicit result than $L^{2}$-local minimality.  The results for the functional $\eone$ are easily carried over to the functionals $\etwo$ and $\ethree$.

We discuss in Section \ref{last} the effect of modifying the surface energy term.  The reasoning set out in Section \ref{last} points out that any path that connects $u=0$ with the global minimizer must, provided $\Delta$ is large enough, pass through a state at which sets $B(u)$ of positive measure first appear.  This is the basis for a calculation which tests whether such states automatically cost a certain minimum amount of energy, analagous to a `nucleation cost'.  The answer is that sets $B$ of positive measure can appear at arbitrarily small energy (measured by any of the $\ei$).  It may help to recall that for a given function $u$ in $\sca_{i}$ the set $B(u)=\{(x,y) \in \om: |v_{y}(x,y)| \geq 1\}$.  We surmise that if there is to be a nucleation cost then it must be as a consequence of some factor beyond the details of the energy functional itself, such as the dynamics, for example.  See Section \ref{last} for further details.    The sets $B(v)$ also play an important role in the local mimimality results of Section \ref{local}.  This is another reason that the introduction of sets $B(u)$ of positive measure is of interest to us.

\subsection{Notation}  One and two dimensional Lebesgue measure are denoted respectively by $\elo$ and $\elt$ throughout the paper.  The usual Sobolev space notation $W^{k,p}$ is used to denote the class of $k-$times weakly differentiable $L^{p}$ functions all of whose derivatives lie in $L^{p}(\om)$.  The $L^{p}$ norm of a function $v$ is denoted by $||v||_{p}$ or $||v||_{L^{p}(\om)}$ depending on the context.  The class of absolutely continuous functions is denoted $AC$, and those functions which are absolutely continuous along almost all lines parallel to the coordinate axes are written $ACL$.  See \cite{Zi} for further details.  The characteristic function of any set $S$ is written $\chi_{S}$. 
All other notation is standard, with the possible exception that the value of the dimensionless positive constants $C$ appearing in various inequalities may, where no confusion arises, change from line to line.  Where it is necessary to distinguish between positive constants we shall use the convention that if $C$ and $c$ appear in the same calculation then $c< C$.  We will also employ the convention that roman letters $\tx$ represent vectors in $\R^{2}$.

\section{The effect of a variable well-depth}\label{variablewell}

When $\Delta=0$ it is clear by inspection that $u=0$ is a global minimizer of $\ei$ in $\sca_{i}$ for $i=1,2,3$.  Therefore the following definition makes sense:
\begin{defn} For each $i=1,2,3,$
\[\Delta_{i}=\sup\{\Delta \geq 0: \ u=0 \ \textrm{globally  minimizes} \ \ei \ \textrm{in} \ \sca_{i}\}.\]
\end{defn}
Each $\Delta_{i}$ will be referred to as a critical well-depth.    The ordering
\begin{equation}\label{ordre} \eone(v) \leq \etwo(v) \leq \ethree(v) \end{equation}
for all appropriate $v$, together with the fact that $\eone(0)=\etwo(0)=\ethree(0)$ for fixed $\eps$ and $\Delta$, implies the inclusions
\begin{eqnarray*} \{\Delta: \ \eone(v) \geq \eone(0) \ \forall v \in \sca_{1}\} & \subset & \{\Delta: \ \etwo(v) \geq \etwo(0) \ \forall v \in \sca_{2}\} \\
\{\Delta: \ \etwo(v) \geq \etwo(0) \ \forall v \in \sca_{2}\} & \subset &  \{\Delta: \ \ethree(v) \geq \ethree(0) \ \forall v \in \sca_{3}\},
\end{eqnarray*}
and hence
\begin{equation}\label{deltaorder}\Delta_{1} \leq \Delta_{2} \leq \Delta_{3}.\end{equation}
It will be shown in this section that all the $\Delta_{i}$ scale alike with respect to $\eps$ and $L$ in the sense that there are dimensionless constants $c <  C$ such that
\[ \frac{c \eps}{L} \leq \Delta_{1} \leq \Delta_{2} \leq \Delta_{3} \leq \frac{C\eps}{L}\]
provided $\eps$ is sufficiently small.

To begin with, Proposition \ref{deltapos} below shows that there is a constant $C$ depending only $L$ such that $\Delta_{1} \geq C \eps^{2}$.  This follows relatively simply by using only the surface energy term in the functional $\eone$ and it turns out to be a crude lower bound on $\Delta_{1}$.  A better (i.e., larger) lower bound is obtained in Section \ref{lb} where it shown that
\[\Delta_{1} \geq \frac{c\eps}{L}.\]
The optimality (in terms of its scaling in $\eps$ and $L$) of this lower bound is proved by evaluating the energy of a particular element $v^{\eps}$ of $\sca_{3}$ in section \ref{ub}.  The structure of $v^{\eps}$ is given in section \ref{ub}.   It is not a branching microstructure, unlike the global minimizer obtained when $\Delta_{i} \sim 1$ and corresponding to the model studied by Kohn and M\"{u}ller.    
The main result of this part of the paper is:
\begin{thm}\label{t1} With $\ei$, $\sca_{i}$ and $\Delta_{i}$ as above, there are dimensionless constants $c<C$ and $\eps_{0}>0$ such that
\[ \frac{c\eps}{L} \leq \Delta_{i} \leq  \frac{C\eps}{L} \ \ \textrm{if} \ \ 0 < \eps < \eps_{0}.\]
Moreover, for all $\Delta$ the global minimizer $U_{i}$ of $\ei$ in $\sca_{i}$ exists, and when $\frac{\Delta}{\Delta_{i}}$ is large enough it satisfies $\elt(B(U_{i})) > 0$. 
\end{thm}

The proof of Theorem \ref{t1} rests on Propositions  \ref{prop2} and \ref{propx} in Sections \ref{lb} and \ref{ub} below.
The lower bound is proved with the help of an interpolation inequality; the upper bound corresponds to the statement concerning the optimality (in a scaling sense) of the lower bound.  The proof of the existence of the global minimizer $U_{i}$ of $\ei$ in $\sca_{i}$ is a relatively straightforward application of the direct method of the calculus of variations.  It is given in an appendix for completeness.  The last assertion of Theorem \ref{t1} can be deduced from the comparison function constructed during the proof of Proposition \ref{propx}.


\subsection{A simple lower bound on $\Delta_{1}$}\label{s2}
We recall that
 \begin{equation*} \eone(u)  =  \int_{\om}\epsilon^{2}u_{yy}^{2}+u_{x}^{2}\,\tdx +\Delta \elt(A(u)) \end{equation*}
 where
 \begin{equation*}A(u)  =  \{(x,y) \in \om: \ |u_{y}(x,y)| < 1 \}. \end{equation*}

The following lemma will be used to show that $\Delta_{1}$ is strictly positive; 
it will also play an important role in Section \ref{local}.

\begin{lem}\label{l1} Let $u \in \sca_{1}$ satisfy $\elt(B(u)) > 0$ and define for each $x \in [0,L]$ 
\begin{eqnarray}
l_{x} & = & \{x \} \times [0,1] \\
\label{pi}\Pi(B)& = & \{x \in [0,L]: \  \scl^{1}(l_{x} \cap B(u)) > 0 \}.
\end{eqnarray}

Then
\begin{equation}\label{e2}\frac{\int_{\om} u_{yy}^{2} \, \emph{\tdx}}{\elt (B(u))} \geq  \frac{4}{\tau(1-\tau)}\end{equation}
where
\[\tau:=\frac{\elt (B(u))}{\elo (\Pi(B(u)))}.\]
\end{lem}

\begin{proof} We begin by remarking that the quantity on the right-hand side of \eqref{e2} is finite under the assumptions of the lemma.   Let $B=B(u)$ for brevity.  Clearly $\elt (B) \neq 0$ implies $\tau > 0$.  By Fubini's Theorem we have that $\tau \leq 1$; the membership of $u$ in 
$\sca_{1}$ further implies $\tau < 1$.   If this were not so then for almost all $x \in [0,L]$ the lines $l_{x}$ would satisfy
\[\elo (l_{x} \cap B) = 1.\]
Then, since $u_{y}$ is absolutely continuous along almost all lines $l_{x}$, we must have for almost all $x$ that either
\[u_{y}(x,y) \geq 1 \ \textrm{for a.e.} \ y \in [0,1] \]
or
\[u_{y}(x,y) \leq -1 \ \textrm{for a.e.} \ y \in [0,1], \]
contradicting the periodic boundary conditions imposed along $y=0$ and $y=1$.  In fact, this argument shows that $\elo (l_{x} \cap B) < 1$
for almost all $x$.

Now we proceed with the proof of inequality \eqref{e2}.  Let $x \in \Pi(B)$.  The argument above shows that we may assume there is at least one open set $\mathcal{Y}_{1} \subset l_{x} \cap B$ on which $u_{y}(x, \cdot) \geq 1$ and at least one other open set $\mathcal{Y}_{2} \subset l_{x} \cap B$ on which $u_{y} \leq -1$.   We may suppose that $y_{1}:=\sup \mathcal{Y}_{1} < \inf \mathcal{Y}_{2}=:y_{2}$, so that the intervening set is $[y_{1},y_{2}]$.  It is easy to check that the minimum of the scalar functional 
\[ f \mapsto \int_{y_{1}}^{y_{2}}(f'')^{2}\,dy \]
among $f \in W^{2,2}([y_{1},y_{2}],\R)$ satisfying $f'(y_{1}) \geq 1$ and $f'(y_{2}) \leq -1$ is 
$\frac{4}{y_{2}-y_{1}}$.  The minimization calculation is of the `free endpoint' kind, so that it prescribes optimal values for the differences
$f(y^{*})-f(y_{1})$ and $f(y_{2})-f(y^{*})$ where $y_{1} < y^{*} < y_{2}$ satisfies $f'(y^{*})=0$.  Note that $y^{*}$ exists because $v_{y} \in \textrm{AC}(l_{x})$ by assumption.   The minimizer is $f(y)=\frac{\left(y-\frac{(y_{1}+y_{2})}{2}\right)^{2}}{y_{1}-y_{2}}$, with $y^{*}=\frac{y_{1}+y_{2}}{2}$.  Since this calculation is elementary we omit the proof.  

The result is 
\begin{equation*} \int_{l_{x}} {u_{yy}}^{2}(x,y) \, dy  \geq \frac{4}{y_{2}-y_{1}}.\end{equation*}
But since 
\begin{eqnarray*}y_{2}-y_{1} & \leq &\elo (A \cap l_{x}) \\
\elo (A \cap l_{x}) & = & 1 - \elo (B \cap l_{x}) \end{eqnarray*}
we must have 
\begin{equation*} \int_{l_{x}} {u_{yy}}^{2}(x,y) \, dy  \geq \frac{4}{1-\elo (B \cap l_{x})}.\end{equation*}
Integrating over $x \in \Pi(B)$ and applying Jensen's inequality gives
\[\int_{\om} u_{yy}^{2} \,\tdx \geq \frac{4\elo(\Pi(B))^{2}}{\int_{\Pi(B)}1-\elo (B \cap l_{x})\,dx }.\]
Dividing by $\elt (B)$ and rearranging yields inequality \eqref{e2}.
\end{proof}

\begin{rem}\label{r1}\emph{The global minimizer $v$ of the functional on the left-hand side of \eqref{e2} is such that $\tau=\frac{1}{2}$. Though easy to construct, $v$ can never belong to $\sca_{1}$ because it violates the boundary condition $u=0$ at $x=0$.}\end{rem}

In the rest of the paper it will be useful to have a label for those elements $u$ of $\sca_{i}$ for which $\elt (B(u)) > 0$ holds.
  \begin{defn}$\sca_{i}^{+}:=\{u \in \sca_{i}: \ \elt(B(u)) > 0 \}$. \end{defn}

\begin{prop}\label{deltapos} $\Delta_{1} \geq C_{2}\eps^{2} \geq C_{1}\eps^{2}>0$, where
 \begin{eqnarray} C_{1} & = & \inf\left\{ \frac{\int_{\om} u_{yy}^{2} \,\emph{\tdx}}{\int_{\om}u_{y}^{2}\,\emph{\tdx} }: \ u \in \sca_{1}^{+}\right\} 
 \\ C_{2} & = & \inf\left\{ \frac{\int_{\om} u_{yy}^{2} \,\emph{\tdx}}{\elt (B(u))}: \ u \in \sca_{1}^{+}\right\}
 \end{eqnarray}
 \end{prop}
\begin{proof}
By Chebychev's inequality, $C_{2} \geq C_{1}$.  Therefore we need only prove $\Delta_{1} \geq C_{2} \eps^{2}$.  By Remark \ref{r1} above $C_{2}=16$, but the infimum is not attained.  Now
\begin{eqnarray*} \eone(u)-\eone(0) &  = & \int_{\om} \eps^{2}u_{yy}^{2} +u_{x}^{2}\,\tdx- \Delta \elt (B) \\
                                & \geq & \left (\frac{\int_{\om} \eps^{2} u_{yy}^{2}\,\tdx}{\elt(B)} -\Delta \right) \elt(B) \\
                                & \geq & (16 \eps^{2} - \Delta ) \elt (B). 
\end{eqnarray*}
So if $\Delta \leq 16 \eps^{2}$ then $u=0$ is a global minimizer of $\eone$, and hence $\Delta_{1} \geq 16\eps^{2}$.
\end{proof}

\subsection{A refined lower bound on $\Delta_{1}$ \label{lb}}

In this section we show that there is a dimensionless constant $c$ such that $\Delta_{1} \geq \frac{c\eps}{L}$  for all $\eps$. 
This improves on (i.e. increases) the lower bound obtained in Section \ref{s2}.  The reason for the improvement is essentially that  the term $\int_{\om}u_{x}^{2}\,\tdx$ is brought into play.

We shall make use of the standard interpolation inequality 
\begin{equation}\label{interp1}\frac{1}{\sigma^{2}}\int_{0}^{1}f_{yy}^{2}\,dy+\sigma^{2}\int_{0}^{1}(f-\rho)^{2}\,dy \geq C \int_{0}^{1}f_{y}^{2}\,dy ,\end{equation}
which holds for some $C > 0$, all $\rho$, all non-zero $\sigma$ and all $f \in W^{2,2}([0,1],\R)$.  (See,  e.g.,  \cite[Section 7.12]{GT}.)
Let $u \in \sca_{1}$, fix $x \in \Pi(B(u))$ and take $\rho=0$, $f(y)=u(x,y)$ in \eqref{interp1} above.  Using the inequality
\[\int_{0}^{1}u_{y}^{2}(x,y) \,dy \geq \elo(l_{x} \cap B(u)),\]
integrating over $x \in \Pi(B)$ and using Fubini's Theorem we obtain
\begin{equation}\label{interp2}
\frac{1}{\sigma^{2}}\int_{\om}u_{yy}^{2}\,\tdx+\sigma^{2}\int_{\om}u^{2}\,\tdx \geq C\elt (B(u))
\end{equation}
Minimizing the left-hand side of \eqref{interp2} over non-zero $\sigma$ we see that
\begin{equation}\label{interp3}\left( \int_{\om}u_{yy}^{2}\,\tdx  \right)^{\frac{1}{2}}\left( \int_{\om}u^{2}\,\tdx  \right)^{\frac{1}{2}}
\geq \frac{C\elt(B(u))}{2}.
\end{equation}
Note that the constant $C$ is independent of the dimensions of the domain $\om$.  

We also need the standard Poincar\'{e} inequality 
\begin{equation}\label{poincare}\int_{\om} u_{x}^{2} \,\tdx \geq \frac{C}{L^{2}} \int_{\om}u^{2} \,\tdx ,\end{equation}
which uses the boundary condition $u=0$ along $x=0$.  The constant $C$ is independent of the domain dimensions.

\begin{prop}\label{prop2} There is a dimensionless constant $c>0$ such that $\Delta_{1} \geq  \frac{c\eps}{L}$ for all $\eps> 0$.  In particular, the lower bound on $\Delta_{1}$ stated in Theorem \ref{t1} holds.
\end{prop}

\begin{proof} Let $u \in \sca_{1}$ and set $B=B(u)$.  By definition of $\eone$ and from inequalities \eqref{interp3} and \eqref{poincare} we have
\begin{eqnarray}\nonumber \eone(u)-\eone(0)& = & \left (\frac{\int_{\om}\eps^{2}u_{yy}^{2} \,\tdx + \int_{\om} u_{x}^{2} \,\tdx}{\elt(B)}-\Delta\right) \elt(B) \\
\label{ineq1}& \geq & \left(\frac{C\eps^{2}\elt(B)}{4\int_{\om}u^{2}\,\tdx} + \frac{C\int_{\om}u^{2}\,\tdx}{L^{2}\elt(B)} - \Delta\right)\elt(B).
\end{eqnarray}
Letting 
\[t=\frac{\int_{\om}u^{2} \,\tdx}{\elt(B)}\]
we see that the right-hand side of \eqref{ineq1} above has the form 
\[ \left(\frac{C\eps^{2}}{4t}+\frac{Ct}{L^{2}} - \Delta\right)\elt(B).\]
the term in brackets is minimized when $t=c\eps L$ for some constant $c$.  From this it follows that any $\Delta \leq \frac{c\eps}{L}$ is such that $u=0$ is a global minimizer of $\ed$. Therefore $\Delta_{1} \geq \frac{c\eps}{L}$.
 \end{proof}

In some cases one can do better than Proposition \ref{prop2}.  The following lemma shows that the lower bound on $\Delta_{1}$ obtained above is correct with constant $c=1$ provided condition \eqref{con2} below holds.  This supplementary condition amounts to a strengthening of the boundary condition along $y=0$ and $y=1$; it is satisfied, for example, by all admissible functions having compact support in $\om$.

\begin{lem}Let $v \in \sca_{1}^{+}$ satisfy 
\begin{equation} \label{con2}\int_{l_{x}}(v_{y}v_{x})_{y}\,dy  = 0    \end{equation} 
for almost every $x$ in $[0,L]$.
Then for almost every $x$ in $\Pi(B)$ 
\begin{equation}\label{wopper1} \int_{\om} \eps^{2} v_{yy}^{2}+v_{x}^{2} \,\emph{\tdx} \geq \eps \elo(l_{x} \cap B).
 \end{equation}
From this it follows that $\Delta_{1} \geq \frac{\eps}{L}$.

Furthermore, if $v\ in \sca_{2}$ satisfies
\begin{equation}\label{wopper2}\int_{\om}v_{yx}^{2}\,\emph{\tdx} \geq \frac{1}{\eps^{2}} \left(\Delta - \frac{\eps}{\elo(\Pi(B))}\right)\elt(B)
\end{equation}
then
\[\etwo(v) \geq \etwo(0).\]
\end{lem}

\begin{proof} First fix $x \in \Pi(B)$ for which \eqref{con2} holds.   Then
\begin{eqnarray*}\int_{\om} \eps^{2} v_{yy}^{2}+v_{x}^{2} \textrm{dx} & \geq &  \int_{[0,x] \times [0,1]} \eps^{2} v_{yy}^{2}+v_{x}^{2}\, \textrm{dx} \\
& \geq & 2\eps\left|\int_{\om}v_{yy}v_{x}\,\textrm{dx}\right| \\     
  & = & 2 \eps \left|\int_{\om} (v_{y}v_{x})_{y}-v_{y}v_{yx} \,\textrm{dx}\right| \\
  & = & 2 \eps \left|\int_{0}^{x} \left\{ \int_{l_{x'}}(v_{y}v_{x})_{y} \,dy - \int_{l_{x'}}\left(\frac{1}{2}v_{y}^{2}\right)_{x}\,dy \right\} \,dx'\right| \\
  & = &  \eps \int_{l_{x}}v_{y}^{2}\,dy \\
  & \geq &   \eps \elo(l_{x} \cap B),
\end{eqnarray*}
where we have applied \eqref{con2} and the boundary condition $v(0,y) = 0$ for $0 \leq y \leq 1$ to pass from the fourth to the fifth line.
Integrating both sides of inequality \eqref{wopper1} over $\Pi(B)$, dividing by $\elo(\Pi(B))$ and inserting the resulting expression into the definition of $\eone(v)$ gives
\[\eone(v) - \eone(0) \geq \left(\frac{\eps}{\elo(\Pi(B))}-\Delta\right)\elt(B), \]
from which the inequality $\Delta_{1} \geq \frac{\eps}{L}$ follows easily.

Inserting the integrated from of \eqref{wopper1} into $\etwo(v)$ gives
\[\etwo(v)-\etwo(0) \geq  \eps^{2}\int_{\om} v_{yx}^{2}\,\tdx + \left(\frac{\eps}{\elo(\Pi(B))}-\Delta\right)\elt(B).\]
Therefore \eqref{wopper2} implies $\etwo(v) \geq \etwo(0)$ as required.  
\end{proof}

\begin{rem}\emph{Any $\Delta$ satisfying $\Delta \leq \frac{\eps}{L}$ forces \eqref{wopper2} to hold.  Therefore inequality \eqref{wopper2} provides a short-cut to the proof that $\Delta_{2} \geq \frac{\eps}{L}$ whenever \eqref{con2} is true.}\end{rem}

\subsection{A sharp upper bound on $\Delta_{3}$} \label{ub}

We show in this section that there is a constant $C$ independent of $\om$ and $\eps$  such that $\Delta_{3} \leq \frac{C\eps}{L}$ if $\eps$ is sufficiently small.   The idea of the proof can be explained as follows.  Let us suppose that for each $\eps > 0$ there is an element $v^{\eps}$ of $\sca_{3}$ with the properties that  
\begin{eqnarray}\label{ansatz1} \int_{\om}\eps^{2}|D^{2}v^{\eps}|^{2} +  ({v^{\eps}_{x}})^{2} \tdx & \leq & C_{1}\eps \\
\label{ansatz2} \elt(B(v^{\eps})) & \geq & C_{2} L.
 \end{eqnarray}
The constants $C_{1}$ and $C_{2}$ should not depend on $\eps$ or $L$.  Let $\Delta < \Delta_{3}$.  Then in particular 
\begin{equation*} \ethree(v^{\eps}) - \ethree(0) \geq   0 \end{equation*}
on the one hand; and, using \eqref{ansatz1} and  \eqref{ansatz2} above,  
\begin{equation*}  C_{1}\eps- \Delta C_{2}L \geq \ethree(v^{\eps})-\ethree(0) \end{equation*}
on the other.  Thus $\Delta \leq \frac{C_{1}\eps}{C_{2}L}$ whenever $\Delta < \Delta_{3}$.  Letting $\Delta \to \Delta_{3}$ yields the desired upper bound.  It remains to prove the existence of a map $v^{\eps} \in \sca_{3}$ with the properties \eqref{ansatz1} and  \eqref{ansatz2}.

\begin{prop} \label{propx}  There exists a map $v^{\eps}$ in the class $\sca_{3}$ and dimensionless constants $C_{1}$ and $C_{2}$ such that \eqref{ansatz1} and \eqref{ansatz2} hold.  In particular, there is a dimensionless constant $C$ such that
\[\Delta_{3} \leq  \frac{C\eps}{L}\]
 whenever $\eps$ is sufficiently small, proving the upper bound on $\Delta_{3}$ stated in Theorem \ref{t1}.
\end{prop} 
\begin{proof}
Let $k,l,h > 0$ and define $H: [0,l]  \to \R$ by
\[H(x) = h -kx \]
where $kl = \frac{h}{2}$. 
Define the function $w$ on $[0,l] \times [0,2h]$ by 

\begin{displaymath} w(x,y)= \left\{\begin{array}{ll} \frac{y^{2}}{2H(x)} & \textrm{if} \ 0 \leq x \leq  l, \ 0 < y \leq H(x) \\
                            y-\frac{H(x)}{2} &  \textrm{if} \ 0 \leq x \leq  l, \ H(x) \leq y \leq 2h - H(x) \\
                                                                                              
2h-\frac{3H(x)}{2} - \frac{(y-2h)^{2}}{2H(x)} & \textrm{if} \  0 \leq x \leq l, \ 2h-H(x) \leq y \leq 2h .\end{array}               \right.
\end{displaymath}

Now extend $w$ to $[0,l] \times [0,4h]$ by reflection in the line $y=2h$, namely 
\begin{equation*}w(x,y)=w(x,4h-y) \ \textrm{if} \ 0 \leq x \leq l, \ 2h \leq y \leq 4h.\end{equation*} 
It can be checked that 
\begin{eqnarray}\label{block1}\int_{0}^{4h}\int_{0}^{l} |D^{2}w|^{2}\,\,\tdx  &  = &  \frac{1}{k}\left(c_{1}+c_{2}k^{2}+c_{3}k^{4}\right) \\
\int_{0}^{4h}\int_{0}^{l} w_{x}^{2}\,\tdx & = & c_{3}k^{2}lh, 
\end{eqnarray}
where the $c_{i}$ are positive dimensionless constants whose precise values are not important.  
Suppose $h$ is chosen so that $N:=\frac{1}{4h}$ is a positive integer.  Extend $w$ by periodicity to $[0,l] \times [0,1]$ and label the resulting function $w$ again.  A computation using \eqref{block1} above together with $ kl = \frac{h}{2}$ gives:
\begin{eqnarray}\label{block2}\int_{0}^{1}\int_{0}^{l} |D^{2}w|^{2}\,\tdx  &  = &  
\frac{1}{h^{2}}\left(c_{1}+c_{2}\left(\frac{h}{l}\right)^{2}+c_{3}\left(\frac{h}{l}\right)^{4}\right) \\
\int_{0}^{1}\int_{0}^{l} w_{x}^{2}\,\tdx & = & c_{4}\frac{h^{2}}{l}. 
\end{eqnarray}

Clearly $w$ is not an element of $\sca$ because it doesn't satisfy the boundary condition at $x=0$.  But we can interpolate between $w(0,y)$ and the function $y \mapsto 0$ as follows.  Define $v: [0,l] \times [0,1] \to \R$ by $v(x,y)=\frac{x}{l}w(0,y)$ and compute directly.

\begin{eqnarray}\label{block3}\int_{0}^{1}\int_{0}^{l} |D^{2}v|^{2}\,\tdx  &  = &  
\frac{1}{h}\left( c_{4}\frac{l}{h}+c_{5}\frac{h}{l}\right) \\
\int_{0}^{1}\int_{0}^{l} v_{x}^{2}\,\tdx & = & c_{6}\frac{h^{2}}{l}. 
\end{eqnarray}

By construction $v$ and $w$ depend only on the parameters $h$ and $l$.
The last step is to glue $v$ and $w$ together to give an element $v^{\eps}$ of $\sca$. Define 

\begin{displaymath} v^{\eps}(x,y) = \left \{\begin{array}{l l } v(x,y; h,l)  &  0 \leq x \leq l, \ 0 \leq y \leq 1 \\ 
w(x-l,y; h,l) & l \leq x \leq 2l, \ 0 \leq y \leq 1. \end{array} \right. \end{displaymath}
It is straightforward to check that 
\begin{equation}\label{block5}\elt (B(v^{\eps}))=\frac{l}{8}.\end{equation}
and that 

Finally we use \eqref{block2}-\eqref{block5} to compute
\begin{eqnarray*}\ethree(v^{\eps}) - \ethree(0) &  = &  \frac{\eps^{2}}{h^{2}}\left(c_{1}+c_{2}\left(\frac{h}{l}\right)^{2}+c_{3}\left(\frac{h}{l}\right)^{4}+c_{4}l + c_{5}\frac{h^{2}}{l}\right) \\
& + &  c_{6}\frac{h^{3}}{l} +c_{7}\frac{h^{2}}{l} - \frac{\Delta l}{8}
\end{eqnarray*}

The domain of $v^{\eps}$ is $\om$ provided we choose $l=\frac{L}{2}$.
Choosing $\frac{h^{2}}{L}=c^{2}\eps$, where $c$ is such that $(4c(\eps L)^{\frac{1}{2}})^{-1} \in \mathbb{N}$ and $|c-1|$ is minimized, and inserting into the above gives
\begin{eqnarray}\label{charlesives}\ethree(v^{\eps}) - \ethree(0) & = &  \frac{\eps}{c^{2}L} \left(c_{1}+c_{2}\frac{\eps}{L} + c_{3} \left(\frac{\eps}{L}\right)^{2}+c_{4}L + c_{5} \eps\right) \\ & + & c_{6}\eps^{\frac{3}{2}}L^{\frac{1}{2}} + c_{7}\eps  -\frac{\Delta L}{16}.
\end{eqnarray}
By ignoring the term in $\Delta$ it can immediately be seen that \eqref{ansatz1} is satisfied.   From \eqref{block5} we have $\elt(B(v^{\eps})) = \frac{L}{16}$, so \eqref{ansatz2} holds.  By the reasoning set out in the lines following \eqref{ansatz1} and \eqref{ansatz2} this concludes the proof.  
\end{proof}

\begin{rem}\emph{The proof of the upper bound on $\Delta_{3}$ can be obtained directly from \eqref{charlesives} as follows.  Simply note that the inequality $\ethree(v^{\eps})  - \ethree(0) \geq 0$ holds because $\Delta < \Delta_{3}$ has been assumed.  In view of \eqref{charlesives} this gives  $\Delta \leq \frac{C \eps}{L}$.  Hence $\Delta_{3} \leq  \frac{C \eps}{L}$.  }\end{rem}

\begin{figure}[ht]
       \centering
        \psfragscanon
        \psfrag{a}{{linear interpolation}}
\psfrag{b}{{$y=1$}}
\psfrag{c}{{$x=L$}}
\psfrag{d}{{$\sim \eps^{\frac{1}{2}}$}}
\psfrag{e}{{$x=\frac{L}{2}$}}
\psfrag{f}{{$v^{\eps}_{y}=-1$}}
\psfrag{g}{{$v^{\eps}_{y}=1$}}
\psfrag{h}{{$0$}}
\includegraphics[scale=0.6]{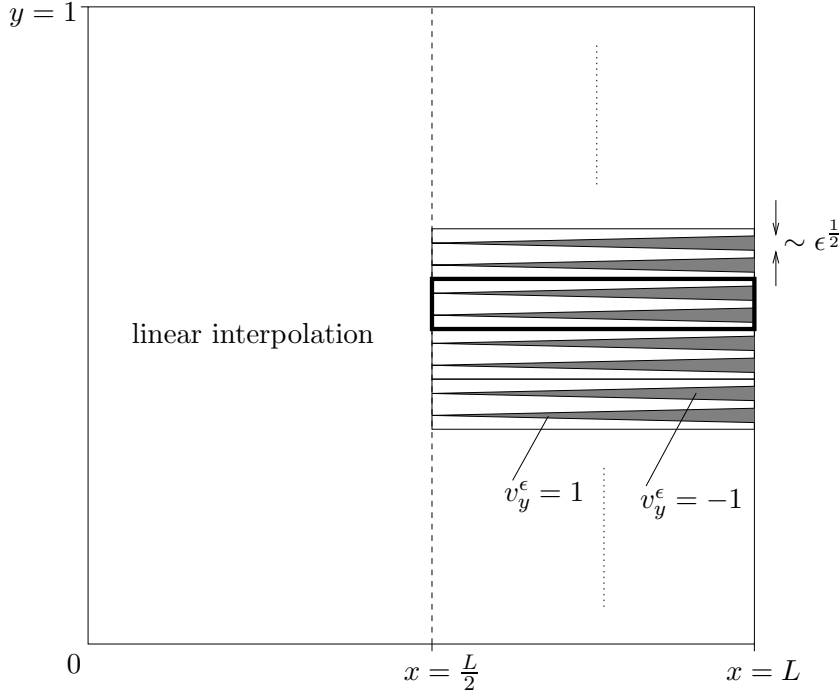}
\caption{The construction of the test function $v^{\eps}$. $|v^{\eps}_{y}|=1$ in the shaded regions; $|v^{\eps}_{y}|<1$ in all other areas.  The basic repeating unit shown with a bold outline in the region $\frac{L}{2} \leq x \leq L$ has a vertical lengthscale of order $\eps^{\frac{1}{2}}$. \label{fig:an1}}
\end{figure}

\section{Austenite as a local minimizer of $\ei$}\label{local}
We saw earlier that $u=0$ is a global minimizer of $\ei$ in $\sca_{i}$ provided $\Delta \leq \Delta_{i}$.  In this section we apply the interpolation inequality \eqref{interp1} to the case $\Delta > \Delta_{i}$ with the aim of proving Theorem \ref{t2} below.  In fact, it suffices to prove the theorem for the functional 
$\eone$ in view of the ordering \eqref{ordre} and since $\eone(0)=\etwo(0)=\ethree(0)$.   The result is then automatically true for the other $\ei$.

\begin{thm}\label{t2}  For each $\eps > 0$ and $\Delta >0$ it is the case that $\eone(v) > \eone(0)$
for nonzero $ v \in \sca_{1}$ such that either
\begin{equation}\label{r} ||v||_{L^{2}(\om)}  <   C \eps^{\frac{7}{2}}\Delta^{-2}\end{equation}
or 
\begin{equation}\label{s} \elt(B(v))  <   C \eps^{6}\Delta^{-4}L^{-1}\end{equation}
holds.
Here, $C$ is a constant independent of $\eps$ and $\Delta$.    The same statements hold with $\ei$ and $\sca_{i}$ in place of $\eone$ and $\sca_{1}$ respectively for $i=2,3$.
\end{thm}  

The spirit of the proof is similar to that of an unpublished result of Ball \cite{Ba06} with the difference that here we take into account the size of the set $B$ where $|u_{y}| \geq 1$.   According to the theorem, $u=0$ is a strict local minimizer in the sense that the  strict inequality $\ei(u) > \ei(0)$ holds whenever $||u||_{L^{2}(\om)}$ or $\elt(B(u))$ is sufficiently small and provided $u \neq 0$.  However, as we shall see in Section \ref{last}, $u=0$ is a degenerate local minimizer of the $\ei$ in the sense that sets $B$ of positive measure can be introduced at arbitrarily small cost measured in terms of $\ei(v)-\ei(0)$ .  It would be interesting to determine whether, in the terminology of Ball \emph{et al} \cite{BKM}, \cite{BM8x},  $u=0$ lies in an energy well of the $\ei$ with respect either to $||v||_{2}$ or $\elt(B(v))$.

We recall from  \eqref{interp1} that for any $v \in \sca$ 
\begin{equation}\label{interp7}
\left(\int_{l_{x}}v_{yy}^{2} \,dy \right)^{\frac{1}{2}} \left(\int_{l_{x}}v^{2} \,dy \right)^{\frac{1}{2}} \geq C \elo(l_{x} \cap B)
\end{equation}
for $a.e. \ x \in [0,L]$ 
where the obvious lower bound $\int_{l_{x}}v_{y}^{2} \,dy \geq \elo(l_{x} \cap B)$ has been used.

We wish to prove $\eone(v) -\eone(0) > 0$ provided $v \neq 0$ and either $\int_{\om}v^{2} \,\tdx$ or $\elt(B(u))$ is small enough.    It is clear that if either
\begin{equation}\int_{\om}v_{yy}^{2}\,\tdx > \frac{\Delta}{\eps^{2}}\elt(B)
 \end{equation}
or \begin{equation} \int_{\om}v_{x}^{2} \,\tdx > \Delta \elt(B)
 \end{equation}
then there is nothing to prove.  There is also nothing to prove should $\elt(B)=0$, since in this case the assumption $v \neq 0$ implies (iin view of the boundary conditions) that $\int_{\om}\eps^{2}v_{yy}^{2}+v_{x}^{2} \,\tdx >0$, and hence that $\eone(v)-\eone(0) > 0$.   Thus we can assume that 
\begin{equation}\label{upper1}\int_{\om}v_{yy}^{2}\,\tdx \leq  \frac{\Delta}{\eps^{2}}\elt(B)
 \end{equation}
and \begin{equation}\label{upper2} \int_{\om}v_{x}^{2} \,\tdx \leq  \Delta \elt(B),
 \end{equation}
where $\elt(B) > 0$.

The claim will be that the two conditions \eqref{upper1} and \eqref{upper2} imply lower bounds on $\int_{\om}v^{2}\,\tdx$ and $\elt(B(u))$, thereby proving the contrapositive of Theorem \ref{t2}. The intuition behind the claim is explained in the course of the next few paragraphs, which should be regarded as a preparation for the proofs of Lemma \ref{l2} and Proposition \ref{l3}. It is on these results that the proof of Theorem \ref{t2} ultimately rests.

\subsection{Preparation for proof of Theorem \ref{t2}} 
Suppose $v \neq 0$ satisfies \eqref{upper1} and \eqref{upper2}.  Applying Lemma \ref{l1} to \eqref{upper1}, and bearing in mind that $\elt(B) >0$, we see that 
\begin{equation}\label{interp5}\frac{4}{t(1-t)} \leq \frac{\Delta}{\eps^{2}},\end{equation}
where 
\[t=\frac{\elt (B)}{\elo (\Pi(B))}.\]
Note that this automatically implies $\Delta \geq 16 \eps^{2}$, which is not a restriction since we already know that $u=0$ is a global minimizer of $\eone$ in $\sca_{1}$ when $\Delta < 16 \eps^{2}$ by Proposition \ref{deltapos}.  Rearranging \eqref{interp5} gives
\begin{equation}\label{proportional1} \frac{4\eps^{2}}{\Delta} \leq \frac{\elt (B)}{\elo (\Pi(B))} \leq 1-\frac{4\eps^{2}}{\Delta},
\end{equation}
where, as usual, $B=B(v)$ for short.  The upper bound merely improves slightly the trivial inequality $t \leq 1$; the lower bound is new information on the set $B(v)$ and is a direct consequence of the assumption \eqref{upper1}.
Inequality \eqref{proportional1} can be interpreted by supposing for the sake of argument that it applies to some rectangle $B$ whose sides are parallel to the coordinate axes. Its `height' would then be bounded below by a fixed constant.  For more general sets $B$ condition eqref{proportional1} should be interpreted in an average sense, \emph{viz.}
\[\frac{1}{\elo(\Pi(B))}\int_{\Pi(B)}\elo(l_{x} \cap B) \,dx \geq \frac{4\eps^{2}}{\Delta}.\]

The following lemma establishes an inequality involving all the terms appearing in the energy $\eone(v)$ and the $L^{2}-$norm of $v$.  It should be regarded as the backbone of Theorem \ref{t2}.

\begin{lem}\label{l2}Let $v \in \sca_{1}^{+}$, let $M>0$ and define 
\begin{eqnarray}\label{piM}\Pi_{M}(B)& = & \left\{x \in \Pi(B): \int_{l_{x}}v_{yy}^{2}(x,y)\,dy  < M\right\} \\
\nonumber B_{M} & =&\{(x,y) \in B: \ x \in \Pi_{M}(B)\}. 
\end{eqnarray}
Then there is a constant $C > 0$ independent of $v$ and the dimensions of $\om$ such that
\begin{equation}\label{killerinterp}
 \left(\int_{\om}v^{2}(x,y)\,\emph{\tdx}\right)^{\frac{1}{2}} 
\left(\int_{\om}{v_{x}}^{2}(x,y)\,\emph{\tdx}\right)^{\frac{1}{2}}
\geq \frac{C}{M}\left(\frac{\elt(B_{M})}{\elo(\Pi_{M}(B))}\right)^{2}
\end{equation}
\end{lem}
\begin{proof}Let $x' \in \Pi_{M}(B)$ and $0 \leq y \leq 1$.  
By applying standard results from the theory of Sobolev functions we may assume without loss of generality that $x \mapsto v^{2}(x,y)$ is weakly differentiable, and hence that 
\begin{equation*}
\int_{0}^{x'}v(x,y)v_{x}(x,y)\,dx = 
\frac{1}{2} v^{2}(x',y) 
\end{equation*}
on using the boundary condition $v(0,y)=0$.
The left-hand side of this inequality is trivially bounded above by 
\begin{equation}\label{interp8}
 U(x',y):=\frac{\sigma^{2}}{2}\int_{0}^{x'}v^{2}(x,y)\,dx + \frac{1}{2\sigma^{2}}\int_{0}^{x'}{v_{x}}^{2}(x,y)\,dx
\end{equation}
for all non-zero $\sigma$; the dependence on $\sigma$ will be minimized out later.
By \eqref{interp7} we have 
\begin{equation*}
 \int_{l_{x'}} v^{2} \,dy \geq \frac{C}{\int_{l_{x'}}{v_{yy}}^{2} \,dy} {\elo(l_{x'} \cap B)^{2}},
\end{equation*}
which in view of the fact that $x' \in \Pi_{M}(B)$ implies
\begin{equation*}
 \int_{l_{x'}} v^{2} \,dy \geq \frac{C}{M} {\elo(l_{x'} \cap B)^{2}}.
\end{equation*}
Integrating both sides of this expression over $\Pi_{M}(B)$ and applying Jensen's inequality to the right-hand side yields
\begin{equation}\label{interp9}
 \int_{\Pi_{M}(B)}\int_{l_{x'}}v^{2}\,dy\,dx' \geq \frac{C}{M}\frac{\elt(B_{M})^{2}}{\elo(\Pi_{M}(B))}.
\end{equation}
The expression on the right is almost the desired lower bound; a factor of 
$\elo(\Pi_{M}(B))$ is missing from the denominator.  But we know that
\begin{equation}\label{interp10} \int_{\Pi_{M}(B)}\int_{0}^{1}U(x',y)\,dy\,dx' \geq 
\int_{\Pi_{M}(B)}\int_{l_{x'}}v^{2}\,dy\,dx',
\end{equation}
and so it remains to estimate the left-hand side of this inequality from above to see that the missing factor can be recovered.
Now 
\begin{eqnarray*}
\int_{\Pi_{M}(B)}\int_{0}^{1}\int_{0}^{x'}v^{2}(x,y)\,\tdx\,dx' & \leq  & 
\int_{\Pi_{M}(B)}\left(\int_{0}^{L}\int_{0}^{1}v^{2}(x,y)\,dy\,dx\right)\,dx' \\
& =  & \elo(\Pi_{M}(B)) \int_{\om}v^{2}\,\textrm{dx}
 \end{eqnarray*}
The term in $u_{x}^{2}$ can be estimated in the same way.  It follows that 
\begin{eqnarray*}
 \elo(\Pi_{M}(B)) \left( 
\frac{\sigma^{2}}{2}\int_{\om}v^{2}(x,y)\,\textrm{dx} +   \frac{1}{2\sigma^{2}}\int_{\om}{v_{x}}^{2}\,\textrm{dx}\right)
& \geq &  \\
\int_{\Pi_{M}(B)}\int_{0}^{1}U(x',y)\,dy\,dx'. & & 
\end{eqnarray*}
Putting this together with \eqref{interp9} and \eqref{interp10} gives
\begin{equation*}
\frac{\sigma^{2}}{2}\int_{\om}v^{2}\,\textrm{dx} +   \frac{1}{2\sigma^{2}}\int_{\om}{v_{x}}^{2}\,\textrm{dx}
\geq \frac{C}{M}\left(\frac{\elt(B_{M})}{\elo(\Pi_{M}(B))}\right)^{2},
\end{equation*}
which on taking $\sigma^{2}= (\int_{\om}v_{x}^{2}\,\textrm{dx})^{\frac{1}{2}}(\int_{\om}v^{2}\,\textrm{dx})^{-\frac{1}{2}}$ concludes the proof.
\end{proof}

We continue to suppose that $v \in \sca_{1}^{+}$ satisfies \eqref{upper1} and \eqref{upper2}. Now by definition of $\Pi_{M}(B)$ it is the case that
\[M \chi_{\Pi(B)\setminus \Pi_{M}(B)}(x) \leq \int_{l_{x}}v_{yy}^{2}\,dy, \]
from which it follows by integrating and then applying \eqref{upper2} that
\begin{equation}\label{interp11}
 \elo(\Pi(B) \setminus \Pi_{M}(B)) \leq \frac{\Delta \eps^{-2}\elt(B)}{M}.
\end{equation}
Hence 
\begin{eqnarray}\nonumber
 \elt(B \setminus B_{M})&  = &  \int_{\Pi_{B} \setminus \Pi_{M}(B)}\elo(l_{x} \cap B)\,dx \\ \nonumber & \leq & \elo(\Pi(B) \setminus \Pi_{M}(B)) \\
\label{interp12}  & \leq & \frac{\Delta \eps^{-2}\elt(B)}{M}.
\end{eqnarray}
We are free to choose $M = 2\Delta \eps^{-2}$, thereby ensuring
\[\elt(B_{M}) \geq \frac{1}{2} \elt(B).\]

Now we combine these observations, an upper bound on $\frac{\int_{\om} v_{x}^{2}\,\textrm{dx}}{\elt(B)}$ and Lemma \ref{l2} to give lower bounds on $\int_{\om}v^{2}\,\textrm{dx}$ and $\elt(B)$.

\begin{prop}\label{l3}Let $v \in \sca_{1}^{+}$ satisfy \eqref{upper1} and \eqref{upper2}. Then there is a constant $C$ independent of $v$ and the dimensions of $\om$ such that
\begin{equation}\label{lower1}
\left(\int_{\om}v^{2}\,\emph{dx}\right)^{\frac{1}{2}} \geq \frac{C\eps^{6}}{(\elt(B))^{\frac{1}{2}}\Delta^{\frac{7}{2}}}.
 \end{equation}
Furthermore, provided $\frac{\Delta}{\Delta_{1}}$ is sufficiently large, 
\begin{eqnarray}\label{lower2}
\left(\int_{\om}v^{2}\,\emph{dx}\right)^{\frac{1}{2}} & \geq &   C \eps^{\frac{7}{2}}\Delta^{-2} \\
\elt(B(v)) & \geq & C \eps^{6}\Delta^{-4}L^{-1}.
\end{eqnarray}

\end{prop}
\begin{proof} Applying Lemma \ref{l2} to $v$ with the choice of $M$ made above and by using inequalities \eqref{upper2} and \eqref{proportional1} we see that 
\begin{eqnarray*}\left(\int_{\om}v^{2}\,\tdx \right)^{\frac{1}{2}} & \geq & \frac{1}{(\int_{\om}{v_{x}}^{2})^\frac{1}{2}}\frac{C}{\Delta \eps^{-2}}\left(\frac{\elt(B)}{\elo(\Pi(B))}\right)^{2} \\
& \geq & \frac{C\eps^{2}}{\Delta} \left(\frac{\eps^{2}}{\Delta}\right)^{2}\frac{1}{(\Delta \elt(B))^{\frac{1}{2}}} \\
& = & \frac{C\eps^{6}}{(\elt(B))^{\frac{1}{2}}\Delta^{\frac{7}{2}}}.
\end{eqnarray*}
This inequality is \eqref{lower1}.  The constant $C$ changes from line to line but it remains independent of $\eps$, $\Delta$ and $L$.

To prove \eqref{lower2} we let $p:=||v||_{L^{2}(\om)}$ and $q:=\elt(B)^{\frac{1}{2}}$ and note that \eqref{lower1} implies 
\begin{equation}\label{pq1}
 pq \geq C \eps^{6}\Delta^{-\frac{7}{2}}.
\end{equation}
For brevity we denote the right-hand side of this inequality by $f$. 

Next, we use the simple interpolation inequality \eqref{interp3} together with \eqref{upper2} to get
\[ C\eps^{2}\Delta \geq   \frac{\elt(B)}{\int_{\om}v^{2}\,\textrm{dx}}.\]  
Hence
\begin{equation}\label{pq2} 
p \geq gq,
 \end{equation}
where $g=C \eps \Delta^{-\frac{1}{2}}$.  
 
Finally, \eqref{poincare} and \eqref{upper2} together imply
\[ \Delta \elt(B) \geq \frac{C}{L^{2}} \int_{\om}v^{2}\,\textrm{dx},\]
which in terms of $p$ and $q$ can be written
\begin{equation}\label{pq3} p \leq C L\Delta^{\frac{1}{2}} q. \end{equation}
The aim is to determine the $(p,q)$ region which is compatible with these inequalities.
This can be done by looking at Fig \ref{fig:pq} 
below.  Note that the line with equation $p=CL\Delta^{\frac{1}{2}}q$ lies above the line with equation $p=gq$ provided $\Delta \geq \frac{C \eps}{ L}$, which, in view of Theorem \ref{t2}, is true whenever $\frac{\Delta}{\Delta_{1}}$ is sufficiently large.  This ensures that the required $(p,q)$ region is nonempty.  

It is immediate that solving $pq=f$ and $p=gq$ yields the smallest possible value $p_{\textrm{min}}$ of $p$ consistent with \eqref{upper1} and \eqref{upper2}.  The result is
\begin{equation*}
p_{\textrm{min}}  =  C \eps^{\frac{7}{2}}\Delta^{-2},
\end{equation*}
giving the lower bound on $||v||_{L^{2}(\om)}$ stated in \eqref{lower2}.  Similarly, the smallest value $q_{\textrm{min}}$ of $q$ consistent with \eqref{upper1} and \eqref{upper2} is found by solving for $q$ in $pq=f$ and $p=CL\Delta^{\frac{1}{2}}q$.  The result is
\begin{equation*}
q_{\textrm{min}}  =  C \eps^{3}\Delta^{-2}L^{-\frac{1}{2}}
\end{equation*}
giving the claimed lower bound on $\elt(B)$.   This concludes the proof of Proposition \ref{l3}.
\end{proof}


\begin{figure}[ht]
        \centering
        \psfragscanon
        \psfrag{labelq}{{$q$}}
\psfrag{labelaone}{{$a_{1}$}}
\psfrag{labelatwo}{{$a_{2}$}}
\psfrag{labelbone}{{$b_{1}$}}
\psfrag{labelbtwo}{{$b_{2}$}}
\psfrag{labellower}{{$p=gq$}}
\psfrag{labelupper}{{$p=CL\Delta^{\frac{1}{2}}q$}}
\psfrag{labelcurve}{{$pq=f$}}
        \psfrag{labelp}{{$p$}}
        \psfrag{labelpmin}{{$p_{\textrm{min}}$}}
        \psfrag{labelqmin}{{$q_{\textrm{min}}$}}
        \includegraphics[scale=0.8]{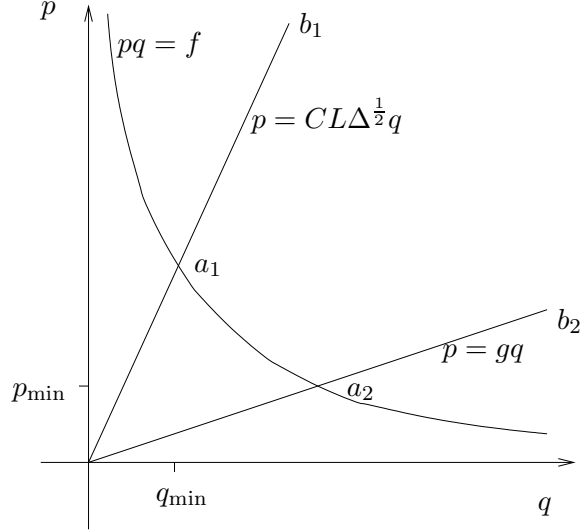}
\caption{The region compatible with \eqref{upper1} and \eqref{upper2} lies within the wedge formed by the lines $a_{1}b_{1}$ and $a_{2}b_{2}$ and above the curve containing the points $a_{1}$ and $a_{2}$.\label{fig:pq}}
\end{figure}

We now draw the preceding results together.

\vspace{2mm}
\noindent{\emph{Proof of Theorem \ref{t2}}   We prove the contrapositive under the assumption $v \neq 0$.  Suppose $\eone(v) \leq \eone(0)$.  Then in particular both
\[\int_{\om} v_{yy}^{2}\,\tdx \leq \frac{\Delta}{\eps} \elt(B) \]
and
\[\int_{\om} v_{x}^{2}\,\tdx \leq \Delta \elt(B)\]
must hold.  These equations are \eqref{upper1} and \eqref{upper2} respectively, where, in view of $v \neq 0$, it can be assumed that $\elt(B) > 0$.  It follows from Proposition \ref{l3} that both $||v||_{L^{2}(\om)} \geq r(\eps, \Delta)$ and $\elt(B) \geq s(\eps, \Delta)$ must hold, concluding the proof of Theorem \ref{t2}.

\begin{rem}\emph{The lines $p=gq$ and $p=CL\Delta^{\frac{1}{2}}q$ coincide when $\Delta \sim \Delta_{1}$ and their relative positions are reversed provided $\frac{\Delta}{\Delta_{1}}$ is small enough.  Under these circumstances the region compatible with all three inequalities \eqref{pq1}, \eqref{pq2} and \eqref{pq3} becomes empty and the starting point for these calculations, namely the inequality 
$\eone(v) \leq \eone(0)$, is contradicted.  But this makes sense since for small $\frac{\Delta}{\Delta_{1}}$ it is the case that $\eone(v) > \eone(0)$ for all non-zero $v$.}\end{rem}

 \begin{rem}\emph{If the lower bounds on $||v||_{L^{2}(\om)}$ and $\elt(B)$ were optimal in a scaling sense then as $\Delta \to \Delta_{1}$ we might expect $r(\eps, \Delta) \to \infty$ and $s(\eps, \Delta) \to \elt(\om)$.  In other words, the constraints on $||v||_{2}$ and $\elt(B(v))$ should become redundant as $\Delta$ approaches $\Delta_{1}$, since when $\Delta \sim \Delta_{1}$ the state $u=0$ is the unique global minimizer and the desired inequality $\eone(v)\geq \eone(0)$ should hold for any admissible $v$.   But it is easily checked that $\Delta \to \Delta_{1}$ implies neither $r(\eps, \Delta) \to \infty$ nor $s(\eps, \Delta) \to \elt(\om)$.  Thus $r(\eps, \Delta)$ and $s(\eps, \Delta)$ would appear to be smaller than they could be, i.e. suboptimal.  Said differently, the inequalities \eqref{r} and \eqref{s} are likely to be sufficient but not necessary conditions for the inequality $\eone(v) \geq \eone(0)$.}
\end{rem}

\section{The effect of modifying the surface energy term}\label{last}


Assuming for argument's sake that $\frac{\Delta}{\Delta_{1}}$ is large enough to ensure, with the aid of Theorem \ref{t1}, that $u=0$ is not the global minimizer of $\ei$ for $i=1,2,3$ then it is clear that, regardless of the dynamics governing the passage from $u=0$ to the global minimizer, 
the appearance of sets $B$ of positive measure is inevitable.   Once such sets have been introduced it ought to be possible to lower the energy by enlarging them, or by allowing them to evolve in some way in order that the energy is driven to its global infimum.   We wish to understand the mechanism behind the introduction of sets $B$ of positive measure into $\om$ using elements of $\sca_{i}$ for $i=1,2,3$.   We suppose that such a mechanism favours small energy.  

More specifically, we are interested in the least value of $\ei(v) - \ei(0)$ consistent with the appearance of sets $B$ of positive measure.   It is clear that unless one restricts the value of $\elt(B)$ the infimum of the difference
$\ei(v)-\ei(0)$ over $\sca_{i}$  will merely reflect the energy of the global minimizer of $\ei$ in $\sca_{i}$.  The results of Section \ref{local} on local minimality imply that if $\elt(B)$ is sufficiently small than $\ei(v) > \ei(0)$.  In particular, therefore, the quantity 
\begin{equation}\label{nucleationcost}
 \liminf_{\mu \to 0} \inf\{\ei(v): v \in \sca_{i}^{+} \ \textrm{such that} \ \elt(B(v)) < \mu\}
\end{equation}
should be nonnegative, and if it were strictly positive then we could interpret it as a lower bound on the cost of `nucleation'. This is investigated below.  It turns out that these costs are in fact zero, even in the case of $\ethree$ where the surface energy includes the full second gradient of $v$.    According to the remarks in the introduction, this may not be because the potential $W$ is not strongly convex in a neighbourhood of the origin.  The force of these remarks was that the potential $W$ can be bounded below by a strongly convex function in a neighbourhood of zero via a simple argument involving a Poincar\'{e} inequality. Instead, it could be that the only way to ensure that the
quantity appearing in \eqref{nucleationcost} is strictly positive is to replace $W_{\Delta}$ with a strongly convex potential whose second gradient is bounded below independently of $\eps$ and $\Delta$.


It happens that low energy competitors can be constructed directly in the case of the functional $\eone$: see Proposition \ref{degeneracy1} below.   The construction of low energy competitors in the case of $\etwo$ and $\ethree$ is indirect and is given in Section \ref{fullsecond} below.

\subsection{Low energy paths for the functional $\eone$}\label{nonucleationcost}
\begin{prop}\label{degeneracy1} Let $\ed'$ be as above and assume $\frac{\Delta}{\Delta_{1}}$ is large enough to ensure that $u=0$ is not the global minimizer of $\eone$ in $\sca_{1}$.  Then  
\begin{equation}\label{naanbread}\liminf_{\mu \to 0} \inf\{\eone(v): v \in \sca_{1}^{+} \ \textrm{such that} \ \elt(B(v)) < \mu\} = 0. 
\end{equation} In other words, sets $B(v)$ of positive measure can be introduced at arbitrarily small energies, as measured by $\eone$. 
\end{prop}
\begin{proof} 
The proof of \eqref{naanbread} is by construction as follows.  Let $a, \delta > 0$, $\lambda >1$, $f(x)=\lambda\left(\frac{x-L}{\delta}\right)^{2}$ and define $v: [L-\delta, L] \times [0,2a] \to \R$ by 
\begin{displaymath} v(x,y)=\left\{\begin{array}{ll} \frac{f(x)y^{2}}{2a} & \textrm{if} \ L-\delta \leq x \leq L, \ 0 \leq y \leq a \\
a f(x)-\frac{f(x)(2a-y)^{2}}{2a} & \textrm{if} \ L-\delta \leq x \leq L, \ a \leq y \leq 2a. \end{array} \right.\end{displaymath} 
Extend $v$ to $[L-\delta, L] \times [0,4a]$ by reflection in $y=2a$, i.e. $v(x,y):=v(x,4a-y)$ if $2a \leq y \leq 4a$.  Finally extend $v$ by zero outside $[L-\delta, L] \times [0,4a]$.   It can then be checked that 
\[\eone(v) \leq C \lambda^{2}\left( \frac{\eps^{2}\delta}{a} + \frac{a^{3}}{\delta}\right),\]
where $C$ is a constant independent of $a, \delta, \eps,\lambda$ and $L$.  It is clear that $\elt(B(v))>0$ for each positive $\delta$ and $a$ and for each $\lambda > 1$. In fact 
\[\elt(B(v))=4a\delta (1-\lambda^{-\frac{1}{2}})^{2},\]
 so given $\mu > 0$ we can, by choosing $\lambda - 1$ sufficiently small and positive, ensure that $\mu > \elt(B(v))>0$ independently of the choice of $a$ and $\delta$.
Taking $a=\delta^{\frac{1}{2}}$ we see that
\[\eone(v) \leq 2C \lambda^{2} \delta^{\frac{1}{2}}.\] 
The conclusion of the proposition follows by letting $\delta \to 0$.
\end{proof}

\subsection{Low energy paths for the full second gradient functional $\ethree$\label{fullsecond}}

In this section it will be convenient to work with the functional $\ethree$.  The same results then hold for $\etwo$ by using the relation \eqref{ordre}.  We seek a sequence $v^{(j)}$ of functions in $\sca_{3}$ satisfying
\begin{itemize} \item[(i)] $\elt(B(v^{(j)}))  >  0$ for all $j$,
\item[(ii)] $\elt(B(v^{(j)})) \to 0$ as $j \to \infty$, and 
\item[(iii)] $\ethree(v^{(j)}) - \ethree(0) \to  0 \ \textrm{as} \ j \to \infty$.
\end{itemize}

One way to do this is to take advantage of the fact that if $f \in L^{2}(\om)$ and if $z$ is the Newtonian potential of $f$ (see e.g. 
 \cite[Chapter 4]{GT}), then $\bigtriangleup \!z=f$ (where $\bigtriangleup$ denotes the Laplacian operator, as usual) and 
\begin{equation}\label{laplacianest} \int_{\om} |D^{2}z|^{2} \,\textrm{dx} = \int_{\om} |f|^{2} \,\textrm{dx}. \end{equation}
(This is part of the Calderon-Zygmund theorem, see e.g. \cite[Theorem 9.9]{GT}.)
When $f$ is sufficiently smooth, for example Lipschitz continuous, the representations
\begin{equation}\label{Newtonian} z(\textrm{x})= \int_{\om}\frac{1}{2 \pi}\ln(|\textrm{x}-\textrm{y}|)f(\textrm{y})\,\textrm{dy}\end{equation}
and
\begin{equation}\label{DNewtonian}\nabla z(\tx) = \int_{\om} \frac{1}{2\pi} \frac{\tx - \ty}{|\tx-\ty|^{2}} f(\ty) \,\textrm{dy}\end{equation}
hold.  One now has control of both $||D^{2}z||_{2}$ and the pointwise behaviour of $\nabla z$ via the function $f$.  The next lemma details an appropriate choice of a sequence of functions $f^{(j)}$ whose corresponding potentials $z^{(j)}$ may be used to satisfy (i), (ii) and (iii) above.   In the following we use the convention that $B_{r}(a)$ denotes the open ball of radius $r$ in $\R^{2}$ centred on $a$.

\begin{lem}\label{degeneracy2}Let the sequence of $L^{2}(B_{2}(0),\mathbb{R})$ functions $f^{(j)}$ be defined by 
\begin{displaymath} f^{(j)}(x,y)= \left\{\begin{array}{l l}2^{j}A_{j}\frac{y}{R} & \textrm{if} \ 0 < R \leq 2^{-j} \\
A_{j}y R^{\alpha_{j}} & \textrm{if} \ 2^{-j} \leq R \leq 1 \\
A_{j} \frac{y}{R} &  \textrm{if} \ 1 \leq R \leq 2,\end{array} \right.\end{displaymath}
where \begin{eqnarray*} A_{j} & = & \frac{k}{j}\\ \alpha_{j} & = &  \frac{1}{j} - 2 \\ R^{2}& = & x^{2}+y^{2}. \end{eqnarray*}
Let $\eta(R)$ be a smooth cut-off function with support in $B_{\frac{3}{2}}(0)$ such that $\eta(R)=1$ if $0 \leq R \leq 1$.  Let $z^{(j)}$ be the Newtonian potential of $f^{(j)}$.  Then the constant $k$ can be chosen so that 
\begin{itemize} \item[(a)] $||f^{(j)}||_{2} \to 0$ as $j \to \infty$, and
\item[(b)] $z^{(j)}_{y}(0,0) \geq (L^{2}+1)^{\frac{1}{2}}$ for all sufficiently large $j$.
\end{itemize}
Here, $L$ is a positive constant.
Let $\psi$ be a smooth cut-off function with support in $B_{\frac{3}{2}}(0)$ and which satisfies $\psi(R) = 1$ if $0 < R \leq 1$.  Then the $C^{2}$ functions $\tilde{z}^{(j)}:=\psi(R)z^{(j)}$ have compact support in $B_{2}(0)$ and they satisfy
\begin{itemize}\item[(c)] $\tilde{z}^{(j)}_{y}(0,0) \geq (L^{2}+1)^{\frac{1}{2}}$ for all sufficiently large $j$, and
\item[(d)] $\int_{B_{2}(0)} |D^{2}\tilde{z}|^{2} \,\emph{\tdx} \to 0$  as $j \to \infty$.
\end{itemize}
\end{lem}
\begin{proof} \textbf{(a)} To see (a) we compute $||f^{(j)}||_{2}^{2}$ directly.
\begin{eqnarray*} \int_{B_{2}(0)}(f^{(j)})^{2} \, \tdx & = & \frac{\pi A_{j}^{2}}{2} + \pi A_{j}^{2}\int_{2^{-j}}^{1}R^{3+2 \alpha_{j}} \,dR + \frac{3\pi A_{j}^{2}}{2} \\
& = & 2\pi A_{j}^{2}+ \frac{3\pi j A_{j}^{2}}{16} \\  
& = & \frac{2\pi k^{2}}{j^{2}}\left(1+\frac{3j}{32 }\right).\end{eqnarray*}
Now (a) follows easily (and independently of the choice of the constant $k$).  

\vspace{2mm}
\noindent \textbf{(b)}
To prove (b) first note that each $f^{(j)}$ is Lipschitz continuous, implying in particular that each $z^{(j)}$ is $C^{2}(B_{2}(0),\mathbb{R})$ and that \eqref{DNewtonian} holds with $B_{2}(0)$ in place of $\om$.     Thus
\begin{eqnarray*}2 \pi z^{(j)}_{y}(0,0) & = &  \int_{B_{2}(0)} -\frac{y}{R}f^{(j)}(\tx) \,\tdx \\ 
& = & -2^{j}\pi A_{j}\int_{0}^{2^{-j}}\,dR  - \pi A_{j}\int_{2^{-j}}^{1}R^{\alpha_{j}+1}\,dR - \pi A_{j} \\
& = & -2\pi A_{j} - \frac{\pi k}{2}.\end{eqnarray*} 
Therefore
\[ z^{(j)}_{y}(0,0) = -\frac{k}{4} - A_{j},  \]
which on choosing $k = -6 (L^{2}+1)^{\frac{1}{2}}$, say, and noting that $A_{j} \to 0$ as $j \to \infty$, implies $z^{(j)}_{y}(0,0) \geq (L^{2}+1)^{\frac{1}{2}}$ for all $j$.   This is part (b) of the lemma.

\vspace{2mm}
\noindent \textbf{(c)} Part (c) follows easily from (b) and the definition of $\tilde{z}^{(j)}$ given above.

\vspace{2mm}
\noindent \textbf{(d)}  By noting that
\[\int_{B_{2}(0)} |D^{2}\tilde{z}^{(j)}| \,\tdx \leq C \int_{B_{2}(0)} |z^{(j)}|^{2}+|\nabla z^{(j)}|^{2} + |D^{2}z^{(j)}|^{2} \,\tdx \]
for some constant $C$, it suffices to prove that $\int_{B_{2}(0)} |z^{(j)}|^{2} \,\tdx$ and $\int_{B_{2}(0)} |\nabla z^{(j)}|^{2} \,\tdx$ converge to zero as $j \to \infty$.  The convergence to zero of the term $\int_{B_{2}(0)} |D^{2}z^{(j)}|^{2} \,\tdx$ is guaranteed by \eqref{laplacianest} and part (a) above.   By the representation \eqref{Newtonian}, standard estimates and Fubini's theorem,
\begin{equation} \int_{B_{2}(0)} {z^{(j)}}^{2} \,\tdx  \leq C  \int_{B_{2}(0)} \left\{\int_{B_{2}(0)} (\ln(|\tx - \ty)|)^{2}\,\tdx \right\}f^{(j)}(\ty)^{2}\,\tdy. \end{equation}
Therefore $||z^{(j)}||_{2}^{2} \to 0$ as $j \to \infty$.    To check the convergence of $||\nabla z^{(j)}||_{2}$ to zero, write 
\begin{equation}\label{greenie}\int_{B_{2}(0)} |\nabla z^{(j)}|^{2} \,\tdx = \int_{\partial B_{2}(0)} z^{(j)}\nabla z^{(j)} \cdot d\nu - \int_{B_{2}(0)} z^{j} 
\bigtriangleup \! z^{(j)} \,\tdx.\end{equation}
Next, note that both $z^{(j)}$ and $\nabla z^{(j)}$ converge uniformly to zero on $\partial B_{2}(0)$, which can be verified by using the fact that each $f^{(j)}$ has compact support in $B_{\frac{3}{2}}(0)$ together with the representations \eqref{Newtonian} and \eqref{DNewtonian}.  The second term in \eqref{greenie} can be estimated by using H\"{o}lder's inequality and \eqref{laplacianest} in that order, giving 
\[\int_{B_{2}(0)} |z^{(j)}  \bigtriangleup \! z^{(j)}| \,\tdx \leq C ||f^{(j)}||_{2}^{2}\]
for some generic constant $C$.
Therefore $||\nabla z^{(j)}||_{2} \to 0$, which concludes the proof.\end{proof}

The next result formalises the statement made at the start of this subsection.    The proof can easily be adapted to show that sets $B(v)$ of positive measure can be introduced into $\om$ in such a way that $\etwo(v)-\etwo(0)$ can be made arbitrarily small.   

\begin{thm} Let $\ethree$ and $\sca_{3}$ be as per \eqref{E3} and \eqref{a3} respectively. 
Then there exist sequences $\{v^{(j)}\} \subset \sca_{3}$ such that 
\begin{itemize}\item[(i)] $\elt(B(v^{(j)})) > 0$ for all sufficiently large $j$,
\item[(ii)] $\elt(B(v^{(j)})) \to 0$ as $j \to \infty$, and  
\item[(iii)] $\ethree(v^{(j)}) - \ethree(0) \to 0$ as $j \to \infty$.
\end{itemize}
In other words, sets $B(v)$ of positive measure can be introduced at arbitrarily small energies, as measured by $\ethree$.
\end{thm}

\begin{proof}Let $P$ be the point $(\frac{L}{2}, \frac{1}{2})$ in $\om$ and define the planar affine map $T$ by 
\[T(x) =  \frac{2(x-P)}{(L^{2}+1)^{\frac{1}{2}}}.\]
Then $T(\om) \subset B_{2}(0)$, $\textrm{lip}(T)=\frac{2}{(L^{2}+1)^{\frac{1}{2}}}$ and we can define 
\[v^{(j)}(x)= \tilde{z}^{(j)}(T(x)) \]
for $x \in \om$, where $\tilde{z}^{(j)}$ is as per Lemma \ref{degeneracy2}.  To check that (i) holds is now straightforward.  Indeed, since
\[ v_{y}^{(j)}(x)=\textrm{lip}(T)\tilde{z}_{y}^{(j)}(T(x))\]  
for all $x$ it follows from part (c) of  Lemma \ref{degeneracy2} that 
\[ v_{y}^{(j)}(p) = \textrm{lip}(T)\tilde{z}^{(j)}(0) \geq 2.\]
Since $v^{(j)}$ is $C^{2}$ it follows that $\elt(\{|v_{y}^{(j)}| \geq 1\}) > 0$, which is statement (i) above.  To see statements (ii) and  (iii) note that
\begin{eqnarray*}\ethree(v^{(j)}) & = &  \int_{\om} \eps^{2} |D^{2}v^{(j)}|^{2} + {v^{(j)}_{x}}^{2} \tdx + \Delta \elt(B(v^{(j)})) \\
& \leq & \int_{B_{2}(0)} \eps^{2}(\textrm{lip}(T))^{2}|D^{2}\tilde{z}^{(j)}|^{2} + (\tilde{z}^{(j)}_{x})^{2} \,\tdy + C \int_{\om}(\tilde{z}^{(j)}_{y})^{2} \,\tdy \\
 & \leq & C \int_{B_{2}(0)} |D^{2}\tilde{z}^{(j)}|^{2} + |\nabla \tilde{z}^{(j)}|^{2} \,\tdy,
 \end{eqnarray*}
where we have applied Chebychev's inequality in the second line.    The right-hand side can now be made arbitrarily small by appealing to part (d) of Lemma \ref{degeneracy2}, proving parts (ii) and (iii) of the theorem.
\end{proof}

\begin{rem}\emph{It may be significant that none of the examples constructed above has small support.  It could be that if we require small support (which is physically reasonable) then there may well be an energy barrier associated with the appearance of sets $B(v)$ of positive measure in $\om$.    It is not clear how or why any such condition should be imposed \emph{a priori}.   Indeed, there are many possible constraints on such sets which may be physically reasonable and yet do not enter into these variational models.
Other factors, such as the asymmetry of the boundary conditions, may also have had a role to play.} \end{rem} 

\section{Appendix} This section is included for completeness only.  We show that the global minimzer $U_{i}$, say, of $\ei$ in $\sca_{i}$ exists.  Only the case $i=2$ is considered here: the others follow by analogy.  

Recall that
\[\etwo(v) = \int_{\om} |\nabla(v_{y})|^{2} + v_{x}^{2} \,\tdx + \Delta \elt(A(v))\]
where
\[A(v) = \{(x,y) \in \R^{2}: \ |v_{y}(x,y)| < 1\}.\]
Now $\etwo$ is bounded below by zero, so its infimum in $\sca_{2}$ exists.  Let $\{v^{(j)}\} \subset \sca_{2}$ be a minimizing sequence.   Then, since $\nabla(v_{y}^{(j)})$ is bounded in $L^{2}$, it follows that there is $w$ in $L^{2}$ such that for a subsequence (and after relabeling) 
$v_{y}^{(j)} \rightharpoonup w$ in $W^{1,2}$.  By the Rellich-Kondrachov compactness theorem we can suppose that the sequence $v_{y}^{(j)}$  converges strongly, and is in particular bounded, in $L^{2}$.   Looking again at $\etwo$ it follows that $v_{x}^{(j)}$ is also  bounded in $L^{2}$.  So $|\nabla v^{(j)} |$ is bounded in $L^{2}$, and hence there is some function $U$ in $W^{1,2}(\om, \R)$ such that
\[v^{(j)} \rightharpoonup U \ \textrm{in} \ W^{1,2}.\]
It follows that $w=U_{y}$.    The trace theorems for Sobolev functions now imply that $U \in \sca_{2}$.   (One could use \cite[Section 4.3, Theorem 1]{EG}, or (A.5) in the appendix of \cite{KM94}, for example.)
Finally, the sequential lower semicontinuity of
\[v \mapsto \int_{\om} v_{x}^{2} \,\tdx\] 
and of 
\[v_{y} \mapsto \int_{\om} |\nabla(v_{y})|^{2}\,\tdx\]
with respect to weak convergence in $W^{1,2}$, together with Fatou's lemma, imply that 
\[\liminf_{j \to \infty} \etwo(v^{(j)})  \geq \etwo(U),\]
concluding the proof of the existence of the global minimizer of $\etwo$ in $\sca_{2}$.

\section{Acknowledgement} I would like to thank Prof. John Ball for \cite{Ba06} and Prof. Stefan M\"{u}ller for sparking my interest in this problem.  This work was begun whilst the author was at the Max Planck Institute for Mathematics in the Natural Sciences, Leipzig in 2006 as part of the MULTIMAT network.  It was completed with the support of an RCUK Academic Fellowship at the University of Surrey.

\end{document}